\documentclass[12pt,twoside,a4paper]{amsart}
\usepackage{amssymb}
\usepackage{graphicx}

\date{\today}

\usepackage[latin1]{inputenc}
\usepackage[T1]{fontenc}

\def\End{{\rm End}}

\def\deg{\text{deg}\,}

\def\dbar{\bar\partial}

\def\R{{\mathbb R}}
\def\C{{\mathbb C}}

\def\Hom{{\rm Hom\, }}

\def\Im{{\rm Im\, }}

\def\U{{\mathcal U}}

\def\xzero{X_0}
\def\xkminus{X_{k-1}}
\def\xk{X_k}
\def\xn{X_n}
\def\sgn{\text{sgn}\;}
\def\codim{\text{codim}\,}
\def\ann{\text{ann}\,}
\def\supp{\text{supp}\,}

\def\scalar{\cdot}
\def\1{\mathbf 1}

\def\Z{{\mathbb Z}}
\def\m{{\mathfrak m}}
\def\J{{\mathcal J}}

\def\q{{\mathfrak q}}
\def\p{{\mathfrak p}}

\def\ass{{\text Ass}}
\def\be{\begin{equation}}
\def\ee{\end{equation}}

\newtheorem{thm}{Theorem}[section]

\newtheorem{cor}[thm]{Corollary}
\newtheorem{prop}[thm]{Proposition}

\theoremstyle{definition}

\theoremstyle{remark}

\newtheorem{preremark}{Remark}
\newtheorem{preex}{Example}

\newenvironment{remark}{\begin{preremark}}{\qed\end{preremark}}
\newenvironment{ex}{\begin{preex}}{\qed\end{preex}}

\numberwithin{equation}{section}

\begin{document}

\title[Residue currents from monomial resolutions]{Residue currents constructed from resolutions of monomial ideals}

\date{\today}

\author{Elizabeth Wulcan}

\address{Mathematical Sciences, Chalmers University of Technology and Mathematical Sciences, G\"oteborg University\\SE-412 96 G\"OTEBORG\\SWEDEN}

\email{wulcan@math.chalmers.se}



\begin{abstract}
Given a free resolution of an ideal ~$J$ of holomorphic functions, one can construct a vector-valued residue current ~$R$, whose annihilator is precisely ~$J$. In this paper we compute ~$R$ in case ~$J$ is a monomial ideal and the resolution is a cellular resolution in the sense of Bayer and Sturmfels. A description of ~$R$ is given in terms of the underlying polyhedral cell complex and it is related to irreducible decompositions of ~$J$.
\end{abstract}

\maketitle
\section{Introduction}
The duality principle for residue currents, due to Dickenstein and Sessa, ~\cite{DS}, and Passare, ~\cite{P}, asserts that a complete intersection ideal of holomorphic functions can be represented as the annihilator ideal of a so-called Coleff-Herrera current, ~\cite{CH}. It has been widely used, for example to obtain effective solutions to division problems, ~\cite{BY}, and explicit versions of the Ehrenpreis-Palamodov fundamental principle, ~\cite{BP}, see also ~\cite{BGVY}. In ~\cite{AW} we generalized the duality principle to general ideals of holomorphic functions by constructing, from a free resolution of an ideal ~$J$, a vector-valued residue current $R$, whose annihilator ideal is precisely ~$J$. This was used to extend several results previously known for complete intersections.  Also, these currents have recently been used by Andersson and Samuelsson, ~\cite{AS}, to obtain new results for $\dbar$-equations on singular varieties.

The degree of explicitness of the current ~$R$ of course directly depends on the degree of explicitness of the free resolution. In case ~$J$ is a complete intersection the Koszul complex is exact and the corresponding current is the classical Coleff-Herrera current, compare to ~\cite{PTY} and ~\cite{A}. In general, though, explicit resolutions are hard to find. In this paper we will focus on monomial ideals, for which there has recently been a lot of work done, see for example the book ~\cite{MS} and the references mentioned therein. We compute residue currents associated with so-called cellular resolutions, which were introduced by Bayer and Sturmfels in ~\cite{BaS}, and which can be nicely encoded into polyhedral cell complexes. Our main result, Theorem ~\ref{frost}, is a complete description of the residue current of a so-called generic monomial ideal.

Because of their simplicity and nice combinatorial description monomial ideals serve as a good toy model for illustrating general ideas and results in commutative algebra and algebraic geometry, see ~\cite{T2} for examples, which make them a natural first example to consider. In ~\cite{W} residue currents of Bochner-Martinelli type, in the sense of ~\cite{PTY}, were computed for monomial ideals, and in ~\cite{AW} and ~\cite{W2}, there are presented some explicit computations of residue currents of certain simple monomial ideals that are not complete intersections. 
Also, many results for general ideals can be proved by specializing to monomial ideals. In fact, recall that the existence of Bochner-Martinelli type residue currents as well as the residue currents in ~\cite{AW} is proved by reducing to a monomial situation by resolving singularities.

We start by considering Artinian, that is, zero-dimensional, monomial ideals in Section ~\ref{artinian_case}. Residue currents associated with general monomial ideals are computed essentially by reducing to this simpler case. A priori, the residue current ~$R$ associated with a cellular resolution of an Artinian monomial ideal has one entry ~$R_\tau$ for each $(n-1)$-dimensional face ~$\tau$ of the underlying polyhedral cell complex. The main technical result in this paper, Proposition ~\ref{description}, asserts that each ~$R_\tau$ is a certain nice Coleff-Herrera current:
\[c~\dbar \Big[\frac{1}{z_1^{\alpha_1}}\Big ]\wedge\ldots\wedge\dbar \Big[\frac{1}{z_n^{\alpha_n}}\Big ],\]
where $\alpha=(\alpha_1,\ldots, \alpha_n)$ can be read off from the cell complex and ~$c$ is a constant.
In particular, if $c\neq 0$ the ideal of functions annihilating ~$R_\tau$, $\ann R_\tau$, is $(z_1^{\alpha_1}, \ldots, z_n^{\alpha_n})$.
A monomial ideal of this form, where the generators are powers of variables, is called \emph{irreducible}. One can show that every monomial ideal can be written as a finite intersection of irreducible ideals; this is called an \emph{irreducible decomposition} of the ideal. Note that an irreducible ideal is primary so an irreducible decomposition of a monomial ideal is a refinement of a primary decomposition. Since one has to annihilate all entries ~$R_\tau$ to annihilate ~$R$, $\bigcap \ann R_\tau$ yields an irreducible decomposition of the ideal $\ann R$, which by the duality principle equals ~$J$, and so the (nonvanishing) entries of ~$R$ can be seen to correspond to components in an irreducible decomposition. In particular, the number of nonvanishing entries are bounded from below by the minimal number of components in an irreducible decomposition.

In general, we can not extract enough information from our computations to determine which entries ~$R_\tau$ that are nonvanishing. Still, for ``most'' monomial ideals we can; if the monomial ideal ~$J$ is generic, which means that the exponents in the set of minimal generators satisfy a certain genericity condition (see Section  ~\ref{prelim} for a precise definition), then Theorem ~\ref{main_artin} states  that ~$R_\tau$ is nonvanishing precisely when ~$\tau$ is a facet of the Scarf complex introduced by Bayer, Peeva and Sturmfels, ~\cite{BPS}. In particular, if the underlying cell complex is the Scarf complex, then all entries of ~$R$ are nonvanishing. The cellular resolution so obtained is in fact a minimal resolution of the generic ideal ~$J$. Theorem ~\ref{minimal_resolution} asserts that whenever the cellular resolution is minimal, the corresponding residue current has only nonvanishing entries. Also, the number of entries is equal to the minimal number of components in an irreducible decomposition.

In Section ~\ref{nonart} we extend the results for Artinian monomial ideals to general monomial ideals. The basic idea is to decompose the residue current into simpler parts, which can be computed essentially as in the Artinian case. 
In ~\cite{AW2} it was shown that the residue current $R$ constructed from a free resolution of the ideal $J$ can be naturally decomposed with respect to the set of associated prime ideals of $J$, $\ass J$; 
\begin{equation}\label{pucko}
R= \sum_{\p\in\ass J} R^\p,
\end{equation}
where $R^\p$ has support on the variety $V(\p)$ of $\p$ and has the so-called \emph{standard extension property (SEP)} with respect to $V(\p)$, which basically means that it is determined by what it is generically on $V(\p)$. Moreover, each $\ann R^\p$ is $\p$-primary, and it turns out that to annihilate ~$R$ one has to annihilate all the currents $R^\p$ and so 
\begin{equation}\label{bach}
J=\ann R=\bigcap_{\p\in\ass J} \ann R^\p
\end{equation}
gives a minimal primary decomposition of $J$.
Now, the simpler currents $R^\p$ associated with a monomial ideal $M$ can be computed by reducing to the Artinian case, using ideas from ~\cite{W}.  The result is a vector of certain simple currents that in particular have irreducible annihilator ideals and that correspond to the $\p$-primary components in an irreducible decomposition of $M$. Our main result, Theorem ~\ref{frost} is a complete description of the residue current associated with a generic monomial ideal ~$M$, generalizing Theorem ~\ref{main_artin}. In particular, we get a decomposition of $R$, which is a refinement of \eqref{pucko} and which corresponds to an irreducible decomposition of ~$M$. 

The technical core of this paper is the proof of Proposition ~\ref{description}, which occupies Section ~\ref{proof_of}. It is inspired by ~\cite{W}, where similar results were obtained for currents of Bochner-Martinelli type corresponding to the Koszul complex. When considering general cellular resolutions the computations get more involved though; in particular, they involve finding inverses of all mappings in the resolution. As in ~\cite{W}, the proof amounts to computing currents in a certain toric variety constructed from the generators of the ideal, using ideas from \cite{BGVY} and \cite{PTY}.

\section{Preliminaries and background}\label{prelim}
Let us start by briefly recalling the construction of residue currents in ~\cite{AW}. Consider an arbitrary complex of Hermitian holomorphic vector bundles over a complex manifold ~$\Omega$, 
\begin{equation}\label{bunt}
0\to E_N\stackrel{f_N}{\longrightarrow}\ldots\stackrel{f_3}{\longrightarrow} E_2\stackrel{f_2}{\longrightarrow}
E_1\stackrel{f_1}{\longrightarrow}E_0,
\end{equation}
that is exact outside an analytic variety ~$Z$ of positive codimension, and suppose that the rank of ~$E_0$ is $1$. In $\Omega\setminus Z$, let ~$\sigma_k$ be the minimal inverse of ~$f_k$, with respect to some Hermitian metric, let $\sigma=\sigma_0+\ldots+\sigma_N$, 
$u=\sigma(I-\dbar\sigma)^{-1}=\sigma+\sigma(\dbar\sigma)+\sigma(\dbar\sigma)^2+\ldots$, and let ~$R$ be the analytic continuation of $\dbar|F|^{2\lambda}\wedge u$ to ${\lambda=0}$, where ~$F$ is any tuple of holomorphic functions that vanishes on ~$Z$. It turns out that ~$R$ is a well defined current taking values in $\End(E)$, where $E=\oplus_k E_k$, which has support on ~$Z$, and which in a certain way measures the lack of exactness of the associated complex of locally free sheaves of $\mathcal O$-modules $\mathcal O(E_k)$ of holomorphic sections of ~$E_k$,
\begin{equation}\label{karv}
0\to \mathcal O(E_N)\stackrel{f_N}{\longrightarrow} \cdots\stackrel{f_2}{\longrightarrow} \mathcal O(E_1)\stackrel{f_1}{\longrightarrow}\mathcal O(E_0).
\end{equation}
In particular, if ~$\J$ is the ideal sheaf $\Im(\mathcal O(E_1)\to\mathcal O(E_0))$ and $\varphi\in\mathcal O(E_0)$ fulfills that the ($E$-valued) current $R\varphi=0$, then locally $\varphi\in\J$. 

Moreover, letting ~$R^\ell_k$ denote the component of ~$R$ that takes values in $\Hom(E_\ell, E_k)$ and $R^\ell=\sum_k R^\ell_k$, it turns out that $R^\ell=0$ for $\ell\geq 1$ is equivalent to that ~\eqref{karv} is exact, in other words that it is a resolution of $\mathcal O(E_0)/\J$, see Theorem 3.1 in ~\cite{AW}. We then write $R_k=R_k^0$ without any risk of confusion. In this case, $R\varphi=0$ precisely when $\varphi\in \J$.

Let us continue with the construction of cellular complexes from ~\cite{BaS}. Let ~$S$ be the polynomial ring $\C[z_1,\ldots,z_n]$ and let $\deg m$ denote the multidegree of a monomial ~$m$ in ~$S$. When nothing else is mentioned we will assume that monomials and ideals are in ~$S$. 

Next, a polyhedral cell complex ~$X$ is a finite collection of convex polytopes (in a real vector space ~$\R^d$ for some ~$d$), the \emph{faces} of ~$X$, that fulfills that if $\tau\in X$ and ~$\tau'$ is a face of ~$\tau$ (for the definition of a face of a polytope, see for example ~\cite{Z}), then $\tau'\in X$, and moreover if ~$\tau$ and ~$\tau'$ are in ~$X$, then $\tau\cap \tau'$ is a face of both ~$\tau$ and ~$\tau'$. The dimension of a face ~$\tau$, $\dim \tau$, is defined as the dimension of its affine hull (in $\mathbb R^d$) and the dimension of ~$X$, $\dim X$, is defined as $\max_{\tau\in X} \dim \tau$. Let ~$\xk$ denote the set of faces of ~$X$ of dimension $(k-1)$ ($\xzero$ should be interpreted as $\{\emptyset\}$). Faces of dimension ~$0$ are called \emph{vertices}. We will frequently identify $\tau\in X$ with its set of vertices. Maximal faces (with respect to inclusion) are called \emph{facets}. A face ~$\tau$ is a \emph{simplex} if the number of vertices, ~$|\tau|$, is equal to $\dim \tau+1$. If all faces of ~$X$ are simplices, we say that ~$X$ is a \emph{simplicial complex}. A polyhedral cell complex $X'\subset X$ is said to be a \emph{subcomplex} of $X$. Moreover, we say that ~$X$ is \emph{labeled} if there is monomial $m_i$ in ~$S$ associated to each vertex $i$. An arbitrary face ~$\tau$ of ~$X$ is then labeled by the least common multiple of the labels of the vertices of ~$\tau$, that is $m_\tau=\text{lcm}\{m_{i}|i\in \tau\}$; $m_\emptyset$ should be interpreted as $1$. Let $\alpha_\tau$ denote $\deg (m_\tau)\in \mathbb N^n$. By ~$\mathbb N$ we mean $0, 1, 2, \ldots$. We will sometimes be sloppy and not differ between the faces of a labeled complex and their labels.

Now, let ~$M$ be a monomial ideal in ~$S$ with minimal generators $\{m_1, \ldots, m_r\}$ (recall that the set of minimal generators of a monomial ideal is unique). Throughout this paper $M$ will be supposed to be of this form if nothing else is mentioned. Moreover, let ~$X$ be a polyhedral cell complex with vertices $\{1,\ldots,r\}$ endowed with some orientation and labeled by ~$\{m_i\}$.
We will associate with ~$X$ a graded complex of free $S$-modules: for $k=0,\ldots, \dim X+1$, let ~$A_k$ be the free $S$-module with basis $\{e_\tau\}_{\tau\in \xk}$ 
and let the differential $f_k:A_k \to A_{k-1} $ be defined by 
\begin{equation}\label{sabel}
f_k: e_\tau  \mapsto  
\sum_{\text{facets }\tau'\subset \tau} \sgn(\tau',\tau)~\frac{m_\tau}{m_{\tau'}}~ e_{\tau'},
\end{equation}
where the sign $\sgn(\tau',\tau)$ ($=\pm 1$) comes from the orientation on ~$X$. Note that $m_\tau/m_{\tau'}$ is a monomial. The complex
\[\mathbb F_X:0\longrightarrow A_{\dim X+1}\stackrel{f_{\dim X+1}}{\longrightarrow} \cdots \stackrel{f_2}{\longrightarrow} A_1\stackrel{f_1}{\longrightarrow} A_0\] is the \emph{cellular complex} \emph{supported on} ~$X$, which was introduced in ~\cite{BaS}. It is exact if the labeled complex ~$X$ satisfies a certain acyclicity condition. More precisely, for $\beta\in \mathbb N^n$ let $X_{\preceq \beta}$ denote the subcomplex of ~$X$ consisting of all faces ~$\tau$ for which $\alpha_\tau\leq \beta$ with respect to the usual ordering in $\mathbb Z^n$. Then ~$\mathbb F_X$ is exact if and only if ~$X_{\preceq \beta}$ is acyclic, which means that it is empty or has zero reduced homology, for all $\beta\in\mathbb N^n$, see Proposition 4.5 in ~\cite{MS}. We then say that ~$\mathbb F_X$ is a \emph{cellular resolution} of ~$S/M$. 

In particular, if ~$X$ is the $(r-1)$-simplex this condition is satisfied and we obtain the classical \emph{Taylor resolution}, introduced by Diana Taylor, ~\cite{Tayl}. Note that if ~$M$ is a complete intersection, then the Taylor resolution coincides with the Koszul complex. If ~$X$ is an arbitrary simplicial complex, ~$\mathbb F_X$ is the more general \emph{Taylor complex}, introduced in ~\cite{BPS}. Observe that if ~$X$ is simplicial the orientation comes implicitly from the ordering on the vertices. 

Recall that a graded free resolution $\cdots \longrightarrow A_k \stackrel{f_k}{\longrightarrow} A_{k-1} \longrightarrow \cdots$ is \emph{minimal} if and only if for each ~$k$, ~$f_k$ maps a basis of ~$A_k$ to a minimal set of generators of $\Im f_k$, see for example Corollary 1.5 in ~\cite{E2}.
The Taylor complex ~$\mathbb F_X$ is a minimal resolution if and only if it is exact and for all ~$\tau\in X$, the monomials ~$m_\tau$ and ~$m_{\tau\setminus i}$ are different, see Lemma ~6.4 in ~\cite{MS}.

Now, to put the cellular resolutions into the context of ~\cite{AW}, let us consider the vector bundle complex of the form ~\eqref{bunt}, where ~$E_k$ for $k=0,\ldots, N=\dim X+1$ is a trivial bundle over ~$\mathbb C^n$ of rank ~$|\xk|$, endowed with the trivial metric, and with a global frame $\{e_\tau\}_{\tau\in \xk}$, and where the differential is given by ~\eqref{sabel}. Alternatively, we can regard ~$f_k$ as a section of $E_k^*\otimes E_{k-1}$, that is, 
\[
f_k=
\sum_{\tau\in \xk} \sum_{\text{facets }\tau'\subset \tau} \sgn(\tau',\tau)~\frac{m_\tau}{m_{\tau'}}~ e_\tau^*\otimes e_{\tau'}.
\]
We will say that the corresponding residue current ~$R$ is associated with ~$X$, and we will use ~$R_\tau$ to denote the coefficient of $e_\tau\otimes e_\emptyset^*$. 
It is well known that the induced sheaf complex ~\eqref{karv} is exact if and only if ~$\mathbb F_X$ is. For example it can be seen from the Buchsbaum-Eisenbud theorem, Theorem 20.9 in ~\cite{E}, and residue calculus - the proof of Theorem 3.1 in ~\cite{AW}. 

Observe that the elements in ~$S$ (holomorphic polynomials) can be regarded as holomorphic sections of ~$E_0$. 
In this paper, by the annihilator ideal of a current ~$T$, $\ann T$, we will mean the ideal in ~$S$ which consists of the elements $\varphi\in S$ for which $R\varphi=0$.

For $b=(b_1,\ldots,b_n)\in \mathbb N^n$ we will use the notation ~$\m^b$ for the irreducible ideal $(z_1^{b_1}, \ldots, z_n^{b_n})$. If $M=\cap_{i=1}^q \m^{b^i}$, for some $b^i\in\mathbb N^n$, is an irreducible decomposition of the monomial ideal $M$, such that no intersectand can be omitted the decomposition is said to be \emph{irredundant}, and the ideals ~$\m^{b^i}$ are then called the \emph{irreducible components} of ~$M$. One can prove that each monomial ideal ~$M$ in ~$S$ has a unique irredundant irreducible decomposition. Giving the irreducible components is in a way dual to giving the generators of the ideal (see Chapter 5 on Alexander duality in ~\cite{MS}), and the uniqueness of the irredundant irreducible decomposition corresponds to the uniqueness of the set of minimal generators of a monomial ideal. This duality will be illustrated in Example ~\ref{lexemple}. 

We will be particularly interested in so-called generic monomial ideals. 
A monomial ~$m'\in S$ \emph{strictly divides} another monomial ~$m$ if ~$m'$ divides ~$m/z_i$ for all variables ~$z_i$ dividing ~$m$. 
We say that a monomial ideal ~$M$ is \emph{generic} if whenever two distinct minimal generators ~$m_i$ and ~$m_j$ have the same positive degree in some variable, then there exists a third generator ~$m_k$ that strictly divides the least common multiple of ~$m_i$ and ~$m_j$. In particular ~$M$ is generic if no two generators have the same positive degree in any variable. Almost all monomial ideals are generic in the sense that those which fail to be generic lie on finitely many hyperplanes in the matrix space of exponents, see ~\cite{BPS}.

We will use the notation $\dbar [1/f]$ for the analytic continuation of $\dbar |f|^{2\lambda}/f$ to $\lambda=0$, and analogously by ~$[1/f]$ we will mean $|f|^{2\lambda}/f|_{\lambda=0}$, that is, just the principal value of ~$1/f$. By iterated integration by parts we have that
\begin{equation}\label{salt}
\int_z
\dbar\Big[\frac{1}{z^p}\Big] \wedge \varphi dz = 
\frac{2\pi i}{(p-1)!}\frac{\partial^{p-1}}{\partial z^{p-1}}\varphi(0).
\end{equation}
In particular, the annihilator of ~$\dbar [1/z^p]$ is ~$(z^p)$. 
We will use the fact that 
\begin{equation}\label{envar}
v^\lambda |z|^{2\lambda}\frac{1}{z^a} \bigg |_{\lambda=0}=
\left [\frac{1}{z^a}\right ] 
~~~~~~~~~ \text{ and } ~~~~~~~~~
\dbar (v^\lambda |z|^{2\lambda})\frac{1}{z^a} \bigg |_{\lambda=0}=
\dbar\left [\frac{1}{z^a}\right ],
\end{equation}
if $v=v(z)$ is a strictly positive smooth function; compare to Lemma 2.1 in \cite{A}.

\section{Artinian monomial ideals}\label{artinian_case}
We are now ready to present our results concerning residue currents ~$R$ associated with cellular complexes of Artinian monomial ideals. We are interested in the component $R^0$, which takes values in $\Hom (E_0, E)$. In fact, when \eqref{karv} is exact $R=R^0$. From Proposition ~2.2 in ~\cite{AW} we know that if ~$M$ is Artinian, then ~$R^0=R^0_n$, where ~$R^0_n$ is a $\Hom(E_0, E_n)$-valued current. Thus, a priori we know that ~$R^0$ consists of one entry $R_\tau ~e_\tau \otimes e_\emptyset^*$ for each $\tau\in\xn$. We will suppress the factor ~$e_\emptyset^*$ in the sequel.

\begin{prop}\label{description}
Let ~$M=(m_1,\ldots,m_r)$ be an Artinian monomial ideal, and let ~$R$ be the residue current associated with the polyhedral cell complex ~$X$ with vertices $\{m_1,\ldots,m_r\}$. Then
\begin{equation}\label{djupbla}
R^0=\sum_{\tau\in\xn} R_\tau ~e_\tau,
\end{equation}
where 
\begin{equation}\label{himmelsbla}
R_\tau=c_\tau ~\dbar \Big[\frac{1}{z_1^{\alpha_1}}\Big ]\wedge\ldots\wedge\dbar \Big[\frac{1}{z_n^{\alpha_n}}\Big ].
\end{equation}
Here ~$c_\tau$ is a constant and $(\alpha_1,\ldots,\alpha_n)=\alpha_\tau$. If any of the entries of ~$\alpha_\tau$ is ~$0$, ~\eqref{himmelsbla} should be interpreted as ~$0$.
\end{prop}
The proof of Proposition ~\ref{description} is given in Section ~\ref{proof_of}. 

Observe that the proposition gives a complete description of ~$R^0$ except for the constants ~$c_\tau$. We are particularly interested in whether the ~$c_\tau$ are zero or not. 
Indeed, note that
\[
\ann \dbar \Big[\frac{1}{z_1^{\alpha_1}}\Big ]\wedge\ldots\wedge\dbar \Big[\frac{1}{z_n^{\alpha_n}}\Big ] = \m^{\alpha_\tau},
\]
so that $\ann R_\tau=\m^{\alpha_\tau}$ if ~$c_\tau\neq 0$. Note in particular that ~$\ann R_\tau$ depends only on ~$c_\tau$ and ~$m_\tau$ and not on the particular vertices of ~$\tau$ nor the remaining faces in ~$X$. 
Furthermore, to annihilate ~$R^0$ one has to annihilate each entry $R_\tau$ and therefore 
\[
\ann R^0=\bigcap_{\tau\in X; ~c_\tau\neq 0} \m^{\alpha_\tau}. 
\]

Now, suppose that the cellular complex ~$\mathbb F_X$ is exact. Then, ~$R=R^0$, and from Theorems 3.1 and 7.2 in ~\cite{AW} we know that 
\[
\ann R=M.
\]
Thus a necessary condition for ~$c_\tau$ to be nonvanishing is that ~$M\subset \m^{\alpha_\tau}$. In general though, Proposition ~\ref{description} does not give enough information to give a sufficient condition, as will be illustrated in  
Example ~\ref{counterscarf}. Below we will discuss two situations, however, in which we can determine exactly which ~$c_\tau$ that are nonzero.

First we will consider generic monomial ideals. For this purpose, let us introduce the  \emph{Scarf complex} ~$\Delta_{M}$ of ~$M$, which is the collection of subsets $I\subset\{1,\ldots,r\}$ whose corresponding least common multiple ~$m_I$ is unique, that is,
\[
\Delta_{M}=\{I\subset\{1,\ldots,r\}|m_I=m_{I'}\Rightarrow I=I'\}.
\]
One can prove that the Scarf complex is a simplicial complex, and that its dimension is a most ~$n-1$. In fact, when ~$M$ is Artinian, ~$\Delta_{M}$ is a regular triangulation of ~$(n-1)$-simplex.
For details, see for example ~\cite{MS}. In ~\cite{BPS} it was proved that if ~$M$ is generic, then the cellular complex supported on ~$\Delta_M$ gives a resolution of ~$S/M$, which is moreover minimal. Furthermore, if ~$M$ in addition is Artinian, then 
\begin{equation}\label{irred}
M=\bigcap_{\tau \text{ facet of } \Delta_M} \m^{\alpha_\tau},
\end{equation}
yields the unique irredundant irreducible decomposition of ~$M$. To be precise, originally in ~\cite{BPS}, a less inclusive definition of generic ideals was used, but the results above were extended in ~\cite{MSY} to the more general definition of generic ideals we use. 

We can now deduce the following.
\begin{prop}\label{coefficients}
Let ~$M\subset S$ be an Artinian generic monomial ideal and let ~$R$ be the residue current associated with the polyhedral cell complex ~$X$. Suppose that ~$\mathbb F_X$ is a resolution of $S/M$. Then ~$c_\tau$ in ~\eqref{himmelsbla} is non-zero if and only if ~$\tau\in\xn$ is a facet of the Scarf complex ~$\Delta_M$.
\end{prop}

\begin{proof}
Suppose that ~$\tau\in \xn$ is not a facet of ~$\Delta_M$. 
We show that $M\not\subset\m^{\alpha_\tau}$, which forces ~$c_\tau$ to be zero.

Let ~$J$ be the largest subset of $\{1,\ldots, r\}$ such that $m_J=m_\tau$. Then for some $j\in J$ it holds that $m_{J\setminus j}=m_\tau$, as follows from the definition of ~$\Delta_M$. If ~$m_j$ strictly divides ~$m_\tau$ then clearly $m_j\notin \m^{\alpha_\tau}$ and we are done. Otherwise, it must hold for some $k\in J\setminus j$ that ~$m_k$ and ~$m_j$ have the same positive degree in one of the variables.
Then, since ~$M$ is generic, there is a generator ~$m_\ell$ that strictly divides the least common multiple of ~$m_j$ and ~$m_k$ and consequently also strictly divides ~$m_\tau$. Hence $m_\ell\notin \m^{\alpha_\tau}$.

On the other hand, since ~\eqref{irred} is irredundant, ~$c_\tau$ has to be nonzero whenever ~$\tau$ is a facet of ~$\Delta_M$.
\end{proof}
Thus, to sum up, Propositions ~\ref{description} and ~\ref{coefficients} yield the following description of the residue current of a generic monomial ideal.
\begin{thm}\label{main_artin}
Let ~$M\subset S$ be an Artinian generic monomial ideal and let ~$R$ be the residue current associated with the polyhedral cell complex ~$X$. Suppose that ~$\mathbb F_X$ is a resolution of $S/M$.
Then
\[
R=\sum_{\tau \text{ facet of } \Delta_M} R_\tau ~ e_\tau,
\]
where ~$\Delta_M$ is the Scarf complex of ~$M$, ~$R_\tau$ is given by ~\eqref{himmelsbla}, and the constant ~$c_\tau$ there is nonvanishing.
\end{thm}

In particular if we choose ~$X$ as the Scarf complex ~$\Delta_M$ we get that all coefficients ~$c_\tau$ are nonzero.

\begin{remark}
Observe that it follows from Theorem ~\ref{main_artin} that ~$X$ must contain the Scarf complex as a subcomplex. Compare to Proposition ~6.12 in ~\cite{MS}.
\end{remark}

An immediate consequence is the following.
\begin{cor}\label{dieu}
Let ~$M\subset S$ be an Artinian generic monomial ideal and let ~$R$ be the residue current associated with the polyhedral cell complex ~$X$. Suppose that ~$\mathbb F_X$ is a resolution of $S/M$. Then
\[
M=\bigcap_{\tau \in X} \ann R_\tau
\] 
is the irredundant irreducible decomposition of ~$M$.
\end{cor}

Another situation in which we can determine the set of nonvanishing constants ~$c_\tau$ is when ~$\mathbb F_X$ is a minimal resolution of ~$S/M$.
Indeed, in ~\cite{M} (Theorem 5.12, see also Theorem 5.42 in ~\cite{MS}) was proved a generalization of \eqref{irred}; if ~$M$ is Artinian and ~$\mathbb F_X$ is a minimal resolution of ~$S/M$, then the irredundant irreducible decomposition is given by
\begin{equation}\label{irred2}
M=\bigcap_{\tau \text{ facet of } X} \m^{\alpha_\tau}.
\end{equation}
Hence, from ~\eqref{irred2} and Proposition ~\ref{description} we conclude that in this case all ~$c_\tau$ are nonvanishing.
\begin{thm}\label{minimal_resolution}
Let ~$M\subset S$ be an Artinian generic monomial ideal and let ~$R$ be the residue current associated with the polyhedral cell complex ~$X$. Suppose that ~$\mathbb F_X$ is a minimal resolution of ~$S/M$. Then 
\[
R=\sum_{\tau \text{ facet of } X} R_\tau ~ e_\tau,
\]
where ~$R_\tau$ is given by ~\eqref{himmelsbla} and the constant 
 ~$c_\tau$ there is nonvanishing. 
\end{thm}
Finally, we should remark, that even though we can not determine the set of non-vanishing entries of the residue current associated with an arbitrary cell complex, we can still estimate the number of nonvanishing entries from below by the number of irreducible components of the corresponding ideal.

Let us now illustrate our results by some examples. First observe that the ideal ~$(z^A)=(z^a=z_1^{a_1}\cdots z_n^{a_n}| a\in A\subset \mathbb N^n)$ in ~$S$ is precisely the set of functions that have support in $\bigcup_{a\in A}(a+\mathbb R^n_+)$, where
\[\supp \sum_{a\in \Z^n} c_a z^a=\{a\in \Z^n | c_a\neq 0\},\] 
and thus we can represent the ideal by this set, see Figure ~\ref{stair}. Such pictures of monomial ideals are usually referred to as \emph{staircase diagrams}. The generators ~$\{z^a\}$ should be identified as the ``inner corners'' of the staircase, whereas the ``outer corners'' correspond to the exponents in the irredundant irreducible decomposition.

\begin{ex}\label{lexemple}
Let us consider the case when ~$n=2$. Note that then all monomial ideals are generic. If ~$M$ is an Artinian monomial ideal, we can write
\[
M=(w^{b_1}, z^{a_2}w^{b_2},\ldots,z^{a_{r-1}}w^{a_{r-1}}, z^{a_r}), 
\]
for some integers 
$a_2<\ldots<a_r$ and $b_1>\ldots>b_{r-1}$.
\begin{figure}\label{stair}
\begin{center}
\includegraphics{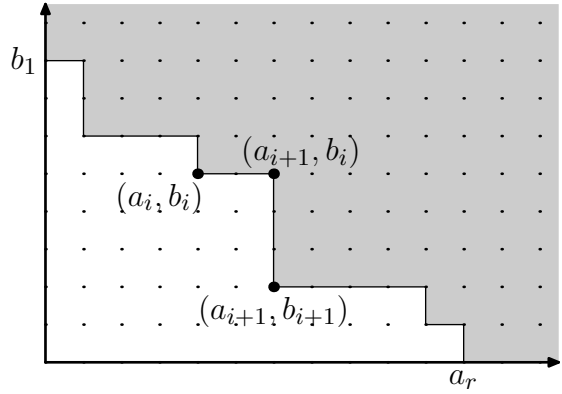}
\caption{The staircase diagram of $M$ in Example \ref{lexemple}.}
\end{center}
\end{figure}
Now ~$\Delta_M$ is one-dimensional and its facets are the pairs of adjacent generators in the staircase. Moreover $m_{\{i,i+1\}}=z^{a_{i+1}} w^{b_i}$, which corresponds precisely to the ~$i$th outer corner of the staircase. 
Thus, according to Theorem ~\ref{main_artin} the residue current ~$R$ associated with a cellular resolution of ~$M$ is of the form 
\[
R=\sum_{i=1}^{r-1}c_i~\dbar\left [\frac{1}{z^{a_{i+1}}}\right]\wedge \dbar 
\left [\frac{1}{w^{b_{i}}}\right] e_{\{i,i+1\}},
\]
for some nonvanishing constants ~$c_i$. The annihilator of the ~$i$th entry is the irreducible component $(z^{a_{i+1}},z^{b_{i}})$. 

Figure ~\ref{stair} illustrates the two dual ways of thinking of ~$M$, either as a staircase with inner corners $(a_i,b_i)$, corresponding to the generators, or as a staircase with outer corners $(a_{i+1},b_i)$, corresponding to the irreducible components or equivalently the annihilators of the entries of ~$R$.
\end{ex}

Let us also give an example that illustrates how we in general fail to determine the set of nonzero ~$c_\tau$ when the ideal is not generic.
\begin{ex}\label{counterscarf}
Consider the non-generic ideal $M=(x^2, xy, y^2, yz, z^2)=:(m_1,\ldots ,m_5)$. The Scarf complex ~$\Delta_M$, depicted in Figure ~\ref{scarffigur}, consists of the 2-simplex $\{2,3,4\}$ together with the one-dimensional ``handle'' made up from the edges $\{1,2\},\{1,5\}$ and $\{4,5\}$.
Moreover the irredundant irreducible decomposition is given by $M=(x, y^2, z)\cap (x^2, y, z^2)$.

Let $X$ be the full 4-simplex with vertices $\{1,\ldots ,5\}$ corresponding to the Taylor resolution. It is then easily checked that for the associated residue current, $c_{\{2,3,4\}}$ and at least one of $c_{\{1,2,5\}}$ and $c_{\{1,4,5\}}$ have to be nonzero, whereas $c_{\{1,2,4\}}$ and $c_{\{2,4,5\}}$ can be either zero or nonzero. The remaining ~$c_\tau$ have to be zero since for them $M\not\subset \m^{\alpha_\tau}$. Thus, in general Proposition ~\ref{description} does not provide enough information to determine which of the coefficients ~$c_\tau$ that vanish.
\begin{figure}\label{scarffigur}
\begin{center}
\includegraphics{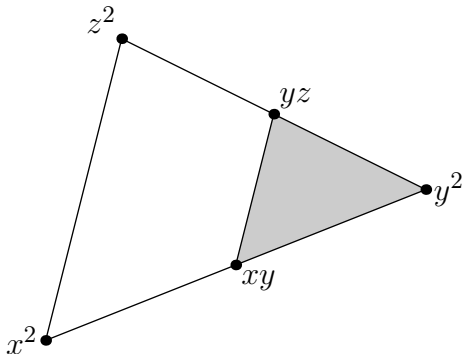}
\caption{The Scarf complex $\Delta_M$ of the ideal $M$ in Example ~\ref{counterscarf}.}
\end{center}
\end{figure}

However, let instead ~$X'$ be the polyhedral cell complex consisting of the two facets $\{2,3,4\}$ and $\{1,2,4,5\}$, that is the triangle and the quadrilateral in Figure ~\ref{scarffigur}. The resolution obtained from ~$X'$, which is in fact the so-called Hull resolution introduced in ~\cite{BaS}, is minimal. Thus, according to Theorem ~\ref{minimal_resolution} the two entries of the associated residue current, which correspond to the two facets of ~$X'$ are both nonvanishing, with annihilators $(x, y^2, z)$ and $(x^2, y, z^2)$ respectively. This could of course be seen directly since we already knew the irredundant irreducible decomposition of ~$M$.
\end{ex}

\section{Proof of Proposition ~\ref{description}}\label{proof_of}
The proof of Proposition ~\ref{description} is inspired by the proof of Theorem ~3.1 in ~\cite{W}, which in turn is inspired by ~\cite{BGVY} and ~\cite{PTY}. We will compute ~$R^0$ as the push-forward of a current on a certain toric manifold constructed from the cell complex ~$X$.

Let us start by giving a description of the current ~$R^0$ in terms of ~$X$. Recall from Section ~2 in ~\cite{AW} that ~$R^0_n$ is the analytic continuation to ~$\lambda=0$ of $\dbar|F|^{2\lambda} \wedge u_n^0$, 
where ~$F$ is a holomorphic function that vanishes at the origin and 
\[u_n^0=(\dbar \sigma_n)(\dbar \sigma_{n-1})\cdots (\dbar \sigma_2) \sigma_1.\] 
Here 
\begin{equation}\label{sigma}
\sigma_k=\frac{\delta_{f_k}^{q_k-1}S_k}{|F_k|^2},
\end{equation}
where ~$q_k$ is the rank of ~$f_k$, $\delta_{f_k}$ is contraction with ~$f_k$, $F_k=(f_k)^{q_k}/q_k!$ and $S_k=(s_k)^{q_k}/q_k!$ is the dual section of ~$F_k$. For details, we refer to Section 2 in ~\cite{AW}, see also ~\cite{W2}.
Furthermore, ~$s_k$ is the section of $E_k\otimes E_{k-1}^*$ that is dual to ~$f_k$ with respect to the trivial metric, that is, 
\[
s_k =\sum_{\tau'\in \xk}~~\sum_{\text{facets }\tau''\subset \tau'} \sgn(\tau'',\tau')~\frac{\overline{m_{\tau'}}}{\overline{m_{\tau''}}} ~e_{\tau'}\otimes e_{\tau''}^*.
\]
Here ~$\overline{m_{\tau'}}$ just denotes the conjugate of ~$m_{\tau'}$. 
Notice that, since $\sigma_k \sigma_{k-1}=0$, as follows by definition, it holds that only the terms obtained when the ~$\dbar$ fall in the numerator survive, and so
\[
u_n^0=
\frac{\dbar(\delta_{f_n}^{q_n-1}S_n)\cdots\dbar(\delta_{f_2}^{q_2-1}S_2)\delta_{f_1}^{q_1-1}S_1}
{|F_n|^2\cdots|F_1|^2}.
\]
Observe furthermore that the numerator of the right hand side of ~\eqref{sigma} is a sum of terms of the form 
\begin{equation}\label{chalmers}
v_k=\pm |\omega_k|^2 \frac{\overline{m_{\tau'}}}{\overline{m_{\tau''}}}~e_{\tau'}\otimes e_{\tau''}^*,
\end{equation}
where $\tau'\in \xk$ and $\tau''\in\xk $ is a facet of ~$\tau'$ and
\[
\omega_k=
\frac{m_{\tau'_1}\cdots m_{\tau'_{q_k-1}}}{m_{\tau''_1}\cdots m_{\tau''_{q_k-1}}},
\]
where for $1\leq \ell\leq q_k-1$, $\tau'_\ell\in \xk$ and $\tau''_\ell\in\xkminus$ is a facet of $\tau'_\ell$. The $\pm$ in front of $|\omega_k|$ depends on the orientation on $X$. Note that the coefficients are monomials. 
It follows that ~$u_n^0$ is a sum of terms of the form 
\[
u_v=u_{\{v_1,\ldots,v_n\}}=\frac{(\dbar v_n)\cdots (\dbar v_2)v_1}{|F_n|^2\cdots |F_1|^2},
\]
where each ~$v_k$ is of the form ~\eqref{chalmers}, and where 
\begin{equation}\label{stockholm}
v_n\cdots v_1=\pm |\omega_n\cdots\omega_1|^2~\overline{m_\tau}~e_\tau\otimes e_\emptyset^*
\end{equation}
for some $\tau\in \xn$.

Observe that each ~$F_k$ has monomial entries. By ideas originally from ~\cite{K} and ~\cite{V}, one can show that there exists a toric variety ~$\mathcal X$ and a proper map $\widetilde\Pi:\mathcal X\to \mathbb C^n$ that is biholomorphic from $\mathcal X\setminus \widetilde\Pi^{-1}(\{z_1\cdots z_n=0\})$ to $\mathbb C^n\setminus \{z_1\cdots z_n=0\}$, such that locally, in a coordinate chart $\U$ of $\mathcal X$, it holds for all ~$k$ that the pullback of one of the entries of ~$F_k$ divides the pullbacks of all entries of ~$F_k$. In other words 
we can write
$\widetilde\Pi^* F_k=F_k^0 F_k'$, where ~$F_k^0$ is a monomial and ~$F_k'$ is nonvanishing, and analogously we have $\widetilde\Pi^*F=F^0 F'$.
The construction is based on the so-called Newton polyhedra associated with ~$F_k$ and we refer to ~\cite{BGVY} and the references therein for details. The mapping ~$\widetilde\Pi$ is locally in the chart ~$\U$ given by
\begin{eqnarray*}
\Pi:\mathcal U &\to &\mathbb C^n\\
t &\mapsto & t^P,
\end{eqnarray*}
where 
$P=(\rho_{ij})$ is a matrix with determinant ~$\pm 1$ and
 ~$t^P$ is a shorthand notation for $(t_1^{\rho_{11}}\cdots t_n^{\rho_{n1}}, \ldots, t_1^{\rho_{1n}}\cdots t_n^{\rho_{nn}})$. Hence, the pullback ~$\Pi^*$ transforms the exponents of monomials by the linear mapping ~$P$;
\begin{equation*}
\Pi^* z^a=\Pi ^* z_1^{a_1}\cdots z_n^{a_n}=
t_1^{\rho_{1}\scalar a}\cdots t_n^{\rho_{n}\scalar a}=
t^{P a},
\end{equation*}
where ~$\rho_i$ denotes the ~$i$th row of ~$P$, so that the pullback of a monomial is itself a monomial. 

Now, from Lemma ~2.1 in ~\cite{AW} we know that $F_k^0\Pi^*\sigma_k$ is smooth in ~$\U$ . However,  
\[
F_k^0\Pi^*\sigma_k=\sum_j \frac{\Pi^* v_k^j}{\bar F_k^0|F'_k|^2}=
\sum_{\alpha\in \mathbb N^n}\sum_{\deg \Pi^* v_k^j=\alpha} \frac{\Pi^* v_k^j}{\bar F_k^0|F'_k|^2},
\]
where ~$v_k^j$ are just the various terms ~$v_k$ that appear in the numerator of ~$\sigma_k$. Therefore clearly for each $\alpha\in \mathbb N^n$ the sum 
\[
\sum_{\deg \Pi^*v_k^j=\alpha} \frac{\Pi^* v_k^j}{\bar F_k^0|F'_k|^2},
\] 
which is just equal to $C t^\alpha/(\bar F_k^0|F'_k|^2)$ for some constant ~$C$, has to be smooth and consequently $t^\alpha/(\bar F_k^0|F'_k|^2)$ is smooth or ~$C=0$. Hence, to compute ~$R^0$ we only need to consider terms ~$u_{v}$, where $v=(v_1,\ldots, v_n)$ is such that $\widetilde\Pi^* v_k/(\bar F_k^0|F'_k|^2)$ is smooth on ~$\mathcal X$ for all ~$k$. 
For such a ~$v$ let us define
\begin{equation}\label{nummer}
R^0_v:=\dbar|F|^{2\lambda}\wedge u_v|_{\lambda=0} ~~~~~\text{ and }~~~~~
\widetilde R^0_v:=\widetilde \Pi^*(\dbar|F|^{2\lambda}\wedge u_v)|_{\lambda=0}. 
\end{equation}
From below it follows that ~$R^0_v$ and ~$\widetilde R^0_v$ are well defined (globally defined) currents and moreover that $\widetilde \Pi_*\widetilde R^0_v=R^0_v$. Furthermore, it is clear that $R^0=\sum R^0_v$, where the sum is taken over all $v$. Next, observe that, in view of ~\eqref{chalmers}, the frame element of $u_v$ is $e_\tau\otimes e_\emptyset^*$, where $\tau\in\xn$ is determined by $v_n$. 
Hence ~$R_\tau ~e_\tau$ in ~\eqref{djupbla} will be the sum of currents ~$R^0_v$, where ~$v$ is such that ~$v_n$ contains the frame element ~$e_\tau$. 
Thus, to prove the proposition it suffices to show that ~$R^0_v$ is of the desired form. 

Let us therefore consider ~$\widetilde R^0_v$ in ~$\U$. Observe that 
\begin{equation}\label{diamond}
\widetilde R^0_v=\dbar |F^0 F'|^{2\lambda}\wedge 
\frac{\Pi^*((\dbar v_n)\cdots(\dbar v_2)v_1)}{|F_n^0\cdots F_1^0|^2 \nu(t)}\bigg|_{\lambda=0}, 
\end{equation}
where $\nu(t):=(|F_n'|\cdots |F_1'|)^2$ is nonvanishing. For further reference, note that ~$\nu(t)$ only depends on $|t_1|, \ldots, |t_n|$. Moreover, let us denote $\deg (F_n^0\cdots F_1^0)\in\mathbb N^n$ by $\gamma=(\gamma_1,\ldots, \gamma_n)$ and $\deg(\omega_n\cdots\omega_1)$ by ~$\beta$, and recall that $\deg m_\tau=\alpha_\tau$.
By Leibniz' rule and \eqref{envar}, recalling \eqref{stockholm}, we see that ~\eqref{diamond} 
is equal to a sum of terms of the form a constant times 
\begin{equation}\label{rubin}
\dbar\bigg[\frac{1}{t_i^{\gamma_i-\rho_i\scalar\beta}}\bigg ]\otimes
\big[\prod_{j\neq i}|t_j|^{2(\rho_j\scalar\beta-\gamma_j)}\big ]
\wedge\frac{\bar t_i^{\rho_i\scalar(\alpha_\tau+\beta)-\gamma_i}
\prod_{j\neq i}\bar t_j^{\rho_j\scalar\alpha_\tau-1}}{\nu(t)}
 \widehat{d\bar t_i}~e_\tau\otimes e_\emptyset^*,
\end{equation}
where ~$t_i$ is one of the variables which fulfills that ~$t_i$ divides the monomials ~$F^0$ and $F_n^0\cdots F_1^0$, whereas $t_1\cdots t_{i-1} t_{i+1}\cdots t_n$ divides ~$\Pi^*m_\tau$. In fact, it is not hard to check that, unless the latter requirement is fulfilled, the corresponding contribution will vanish for symmetry reasons. Moreover ~$\widehat{d\bar t_i}$ is just shorthand for 
$d\bar t_1\wedge \ldots \wedge d\bar t_{i-1} \wedge d\bar t_{1+1} \wedge \ldots \wedge d\bar t_n $. 
Note that since $\Pi^*v_k/(F_k^0|F'_k|^2)$ is smooth there will be no occurrences of any of the coordinate functions ~$\bar t_j$ in the denominator, except for them in ~$\nu(t)$, and in particular it follows that $\gamma_j-\rho_j\scalar \beta\geq 0$ when $j\neq i$. 
Moreover, due to \eqref{salt}, ~\eqref{rubin} vanishes whenever there is an occurrence of ~$\bar t_i$ in the numerator. Hence a necessary condition for ~\eqref{rubin} not to vanish is that 
\[
\rho_i\scalar(\alpha_\tau+\beta)-\gamma_i=0.\]

We will now compute the action of ~$\widetilde R^0_v$ on the pullback of a test form $\phi=\varphi(z) dz$ of bidegree $(n,0)$. Here ~$dz=dz_1\wedge\ldots\wedge dz_n$.
Let $\{\mathcal U_\ell\}$ be the cover of ~$\mathcal X$ that naturally comes from the construction of $\mathcal X$ as described in the proof of Theorem 3.1 in ~\cite{W}, and let ~$\{\chi_\ell\}$ be a partition of unity on ~$\mathcal X$ subordinate $\{\mathcal U_\ell\}$. 
It is not hard to see that we can choose the partition in such a way that the ~$\chi_\ell$ are circled, that is, they only depend on $|t_1|, \ldots, |t_n|$. Now $\widetilde R^0_v=\sum_{\ell}\chi_\ell \widetilde R^0_v$. We will start by computing the contribution from our fixed chart ~$\mathcal U$ (with corresponding cutoff function ~$\chi$), where ~$\widetilde R^0_v$ is realized as a sum of terms ~~\eqref{rubin}.

Recall that ~$R$ has support at the origin; hence it only depends on finitely many derivatives of ~$\varphi$ at the origin. Moreover we know that ~$\bar h$ annihilates ~$R$ if ~$h$ is a holomorphic function which vanishes on ~$Z$, see Proposition 2.2 in ~\cite{AW}. For that reason, to determine ~$R^0_v$ it is enough to consider the case when ~$\varphi$ is a holomorphic polynomial. We can write ~$\varphi$ as a finite Taylor expansion,
\begin{equation*}
\varphi=
\sum_{a}
\frac{\varphi_{a}(0)}{a !} z^a, 
\end{equation*}
where $a=(a_1,\ldots,a_n)$, 
$
\varphi_{a}=
\frac{\partial^{a_1}}{\partial z_1^{a_1}}\cdots
\frac{\partial^{a_n}}{\partial z_n^{a_n}}\varphi
$
 and $a!=a_1!\cdots a_n!$,
with pullback to ~$\mathcal U$ given by
\begin{equation*}
\Pi^* \varphi=
\sum_{a}\frac{\varphi_{a}(0)}{a !}
t^{Pa}
=\sum_{a}\frac{\varphi_{a}(0)}{a !}
t_1^{\rho_1\scalar a} \cdots t_n^{\rho_n\scalar a}.
\end{equation*}
Moreover a computation similar to the proof of Lemma 4.2 in ~\cite{W} yields
\begin{equation*}
\Pi^* dz=
\det P ~ t^{(P-I)\1} ~dt,
\end{equation*}
where ~$\1=(1,1,\ldots,1)$.

Since $\det P\neq 0$, it follows that $\chi \widetilde R^0_v.\Pi^*\phi$ is equal to a sum of terms of the form a constant times 
\begin{multline*}
\int
\dbar 
\bigg [\frac{1}{t_i^{\rho_i\scalar\alpha_\tau}}\bigg]
\otimes 
\big[\prod_{j\neq i}|t_j|^{2(\rho_j\scalar\beta-\gamma_j)}\big ]
\wedge
\frac{\prod_{j\neq i}\bar t_j^{\rho_j\scalar\alpha_\tau-1}}{\nu(t)}
\widehat{d\bar t_i} ~e_\tau\otimes e_\emptyset^*~\wedge
\\
\chi(t)
\sum_{a}
\frac{\varphi_{a}(0)}{a !}t^{Pa} 
t^{(P-I)\1} dt=
~\sum_{a} I_{a}~\wedge
\frac{\varphi_{a}(0)}{a !}
~e_\tau\otimes e_\emptyset^*,
\end{multline*}
where
\begin{equation}\label{fredag}
I_{a}
=
\int
\dbar
\bigg [\frac{1}{t_i^{\rho_i\scalar (\alpha_\tau-a-\1)+1}}\bigg]
\otimes
[\mu_{a}]\wedge
\frac{ \chi(t)}{\nu(t)}~
\widehat{d\bar t_i}\wedge dt.
\end{equation}
Here ~$\mu_{a}$ is the Laurent monomial 
\begin{equation*}
\mu_{a}=
\prod_{j\neq i}
t_j^{\rho_j\scalar(\beta+a+\1)-\gamma_j-1}~
\bar t_j^{\rho_j\scalar(\beta+\alpha_\tau)-\gamma_j-1}.
\end{equation*}
Invoking ~\eqref{salt} we evaluate the $t_i$-integral.
Since ~$\nu$ and ~$\chi$ depend on $|t_1|, \ldots, |t_n|$ 
it follows that
$\frac{\partial^{\ell}}{\partial t_i^\ell}\frac{\chi}{\nu(t)}|_{t_i=0}=0$ for $\ell \geq 1$ and thus ~\eqref{fredag} is equal to
\begin{equation}\label{skynda}
2\pi i
\int_{\widehat t_i}
\frac{\chi(t)|_{t_i=0}
[\mu_{a}]}
{\nu(t)|_{t_i=0}}~
\widehat{d\bar t_i}\wedge \widehat{dt_i},
\end{equation}
if 
\begin{equation}\label{fearful}
\rho_i\scalar (\alpha_\tau - a- \1)+1=1,
\end{equation}
and zero otherwise. 
Moreover, for symmetry reasons, ~\eqref{skynda} vanishes unless
\begin{equation}\label{trip}
\rho_j\scalar(\alpha_\tau-a-\1)=0
\end{equation}
for $j\neq i$, that is, unless ~$\mu_a$ is real.

Thus, since ~$P$ is invertible, the system of equations 
~\eqref{fearful} and ~\eqref{trip}
has the unique solution $a=\alpha_\tau-\1$ if $\alpha_\tau\geq\1$. Otherwise there is no solution, since ~$a$ has to be larger than $(0,\ldots,0)$.
With this value of ~$a$ the Laurent monomial ~$\mu_{a}$ is nonsingular and so the integrand of ~\eqref{skynda},
\begin{equation*}
\frac{\chi(t)|_{t_i=0}~\prod_{j\neq i}|t_j|^{2(\rho_j\scalar(\beta+\alpha_\tau)-\gamma_j-1)}}
{\nu(t)|_{t_i=0}},
\end{equation*}
becomes integrable. Hence ~$I_a$ is equal to some finite constant if $a=\alpha_\tau-\1$ and zero otherwise. 

Now, recall that the chart ~$\mathcal U$ was arbitrarily chosen. Thus adding contributions from all charts reveals that ~$R^0_v$ and thus ~$R_\tau$ is of the desired form ~\eqref{himmelsbla}, and so Proposition ~\ref{description} follows.

\begin{remark}
We should compare Proposition ~\ref{description} to Theorem 3.1 in ~\cite{W}. It states that the residue current of Bochner-Martinelli type of an Artinian monomial ideal is a vector with entries of the form ~\eqref{himmelsbla}, but it also tells precisely which of these entries that are non-vanishing. If we had not cared about whether a certain entry was zero or not we could have used the proof of Proposition ~\ref{description} above. Indeed, the Koszul complex, which gives rise to residue currents of Bochner-Martinelli type, can be seen as the cellular complex supported on the full ~$(r-1)$-dimensional simplex with labels $m_\tau=\{\prod_{i\in F} m_i\}$. It is not hard to see that the proof above goes through also with this non-conventional labeling.
\end{remark}

\section{General monomial ideals}\label{nonart}
If the monomial ideal $M\subset S$ is of positive dimension, the computation of the residue current $R$ associated with a cellular resolution of $S/M$ gets more involved. In general $R=R_p+\ldots + R_\mu$, where $R_k$ has bidegree $(0,k)$ and takes values in $\Hom(E_0,E_k)$, $p=\codim M$ and $\mu=\min(n,r)$. 
Our strategy is to decompose $R$ into the simpler currents $R^\p$, compare to \eqref{pucko}, which can be computed essentially as $R=R_n$ in the Artinian case following the proof of Theorem 5.2 in ~\cite{W}.

In ~\cite{AW2} we introduce a class of currents that we call pseudomeromorphic. The definition is modeled on the residue currents that appear in various works such as \cite{A} and \cite{PTY};  basically a current is pseudomeromorphic if it is the direct image of certain simple ``semi-meromorphic'' currents, see ~\cite{AW2} for precise statements. In particular, all currents in this paper are pseudomeromorphic, and moreover, $\dbar$-closed pseudomeromorphic currents correspond to so-called Coleff-Herrera currents or locally residual currents, compare to \cite{Bj} or \cite{DS}. An important property of pseudomeromorphic currents is that they allow for multiplication with characteristic functions of constructible sets. We use the notation $T|_W$ for the current $\mathbf 1_W T$ if $W$ is a constructible set; indeed, multiplication with $\mathbf 1_W$ can be thought of as restricting to $W$. The current $R^\p$ is then defined as $R|_{V(\p)\setminus \bigcup_{\q\in\ass J, \supset\p}V(\q)}$.

As in the Artinian case, we start by presenting a technical proposition. In \cite{AW2} it is shown that $R^\p$ has the so-called Standard Extension Property (SEP), that is, it is equal to its own standard extension in the sense of \cite{Bj}, which basically means that it has no mass concentrated to subvarieties of its support. As a consequence $R^\p$ behaves essentially like its component of lowest degree, that is, $R^\p_\ell$ if $\codim \p=\ell$. Indeed, outside the set $Z_{\ell+1}\cup\cdots\cup Z_N$, where $Z_j$ denotes the set where the mapping $f_j$ in \eqref{bunt} does not have optimal rank, $R^\p=\beta R^\p_\ell$, where
\begin{equation}\label{beta}
\beta := \sum_{j\geq 0} (\dbar \sigma)^j
\end{equation}
with $(\dbar \sigma)^0$ interpreted as the identity map on $E$, is smooth. 
By arguments similar to the proof of Proposition ~2.2 in ~\cite{AW2} one can show that $\beta T$ has a standard extension over $Z_{\ell+1}\cup\cdots\cup Z_N$ for any $E_\ell$-valued pseudomeromorphic current $T$, see also ~\cite{Bj}. By the Buchsbaum-Eisenbud theorem, see ~\cite{E}, $\codim (Z_{\ell+1}\cup\cdots\cup Z_N)\geq \ell+1$ when \eqref{karv} is exact, and so, since $R^\p$ has the SEP with respect to $V(\p)$ of codimension $\ell$, $R^\p$ is equal to the standard extension of $\beta R^\p_\ell$. 

Recall that each associated prime ideal of a monomial ideal is monomial and therefore generated by a subset of the variables, see for example ~\cite{Sw}. For $K\subset \{1,\ldots, n\}$ let $\p_K$ denote the prime ideal $(z_i)_{i\in K}$. If $\p_K$  of codimension $\ell$ is associated with $M$ we have a priori that $R^{\p_K}_\ell$ consists of one entry $R^{\p_K}_{\ell,\tau} e_\tau\otimes e_\emptyset ^*$ for each $\tau\in X_\ell$. It follows from above that $R^{\p_K}=\sum_{\tau\in X_\ell} R_{(K,\tau)} ~e_\tau\otimes e_\emptyset ^*$, where $R_{(K,\tau)}$ is the standard extension of $\beta R^{\p_K}_{\ell,\tau}$. We will suppress factor $e_\emptyset ^*$ in the sequel. 

\begin{prop}\label{projicerat}
Let $M\subset S$ be a monomial ideal and let $R$ be the residue current associated with the polyhedral cell complex $X$. Suppose that $\mathbb F_X$ is a resolution of $S/M$. 
Then
\[
R=\sum_{\p_K\in\ass M} R^{\p_K},
\]
where 
\[
R^{\p_K} = \sum_{\tau\in X_\ell} R_{(K,\tau)}~e_\tau 
\]
if $\p_K=(z_{k_1}, \ldots, z_{k_\ell})$ is of codimension $\ell$. 
Here $R_{(K,\tau)}$ is the standard extension of 
\begin{equation}\label{himmelsgron}
\beta ~ C_{(K,\tau)}(\eta)\otimes \dbar \Big[\frac{1}{z_{k_1}^{\alpha_{k_1}}}\Big ]\wedge\ldots\wedge\dbar \Big[\frac{1}{z_{k_\ell}^{\alpha_{k_\ell}}}\Big ],
\end{equation}
where $\beta$ is given by \eqref{beta}, $\eta$ denotes the variables $z_i, i\notin K$, 
$C_{(K,\tau)}(\eta)$ is a smooth function outside a set of codimension $>\ell$, and $(\alpha_1,\ldots,\alpha_n)=\alpha_\tau$. If any of the entries of $\alpha_{k_1},\ldots,\alpha_{k_\ell}$ is $0$, \eqref{himmelsgron} should be interpreted as $0$.
\end{prop}

\begin{proof}
Throughout this proof we will use the notation from the proof of Proposition ~\ref{description}. We will start by computing $R_n|_{\{0\}}=R_n^0|_{\{0\}}$; note that if the maximal ideal $\m=(z_1,\ldots, z_n)$ is associated with $M$, then this is precisely $R^\m_n=R^\m$. Following the proof of Proposition ~\ref{description}, we have that $R_n^0=\sum R^0_v$, where $R^0_v$ is given by \eqref{nummer}. Furthermore $R^0_v=\widetilde \Pi_* \widetilde R^0_v$, where $\widetilde R^0_v$ is locally a finite sum of currents $T_\ell$ of the form \eqref{rubin}. Now, according the proof of Proposition ~2.2 in ~\cite{AW2}, $R^0_v|_{\{0\}}=\sum \widetilde \Pi_*  T_{\ell'}$ where the sum is taken over $\ell'$ such that $t_i$ in $T_{\ell'}$ divides $\Pi^{-1}(\{0\})$. Moreover, $R^0_v|_{\{0\}}$ is a pseudomeromorphic current that has support at the origin; in particular it follows that it is annihilated by $\bar h$ for all holomorphic functions $h$ that vanish at the origin, see Proposition 2.3 in ~\cite{AW2}. Hence the rest of the proof of Proposition ~\ref{description} goes through and we get that $R_n|_{\{0\}}$ is of the form \eqref{djupbla}, and so $R^\m$ is of the desired form.

Now, consider an associated prime $\p_K$ of codimension $\ell$. We want to compute $R^{\p_K}_\ell$, which is equal to $R_\ell|_{V_K}$, where $V_K=V(\p_K)=\{z_{k_1}=\ldots=z_{k_\ell}=0\}$, since by Corollary ~2.4 in ~\cite{AW2} a pseudomeromorphic current of bidegree $(p,q)$ that has support on a variety of codimension $>q$ vanishes. Let $W_K$ denote the Zariski-open set $\{z_i\neq 0\}_{i\notin K}$. Note that $\codim(V_k\setminus W_k)=\ell+1$. Thus, since $R_\ell$ is a pseudomeromorphic current of bidegree $(0,\ell)$ we have that $R_\ell|_{V_K}=R_\ell|_{V_K}|_{W_K}$.

The current $R_\ell|_{V_K}|_{W_K}$ can be computed analogously to $R_n|_{\{0\}}$ above. As in the proof of Proposition ~\ref{description} we get that $u_\ell^0$ is a sum of terms 
\[
u_v=\frac{(\dbar v_\ell)\cdots (\dbar v_2)v_1}{|F_\ell|^2\cdots |F_1|^2},
\]
where $v_k$ is given by \eqref{chalmers} and $F_k=(f_k)^{q_k}/q_k!$. 
Denote the $z_i, i\in K$ by $\zeta$ and the $z_i, i\notin K$ by $\eta$ and let $\phi$ be a test form of bidegree $(n,n-\ell)$ with compact support in $W_{K}$. Then $R_\ell$ acting on ~$\phi$ is the analytic continuation to $\lambda=0$ of a sum of terms 
\begin{equation}\label{fosfor}
\int_\eta \int_\zeta \dbar|F|^{2\lambda} \wedge u_v \wedge \varphi(\zeta, \eta) d\zeta\wedge d\bar \eta\wedge d\eta,
\end{equation}
where $v$ fulfills that $\Pi^*v_k/(\bar F_k^0|F_k'|^2)$ is smooth, compare to Section ~\ref{proof_of}. It is easily checked that \eqref{fosfor} vanishes unless $\phi$ is of the form $\varphi(\zeta, \eta) d\zeta\wedge d\bar \eta\wedge d\eta$. We can now compute the inner integral of \eqref{fosfor} as we computed $R_n|_{\{0\}}$ above. Indeed, $V_K$ corresponds to the origin in the $\zeta$-plane and since $\eta$ is nonvanishing in $W_K$ we can regard the coefficients of $v_k$ as monomials in $\zeta$ times monomials in the parameters $\eta$.
Thus we get that the inner integral is of the form
\[
C(\eta) \otimes \dbar \Big[\frac{1}{z_{k_1}^{\alpha_{k_1}}}\Big ]\wedge\ldots\wedge\dbar \Big[\frac{1}{z_{k_\ell}^{\alpha_\ell}}\Big ] \wedge e_\tau\wedge \varphi(\zeta, \eta)d\zeta,
\]
where $C$ is smooth. Hence summing over all $u_v$ gives that $R_\ell|_{V_K}|_{W_K}$, and consequently also $R^{\p_K}$, is of the desired form. Note that $C(\eta)$ extends as a distribution over $W_K^C$.  
\end{proof}

Observe that Proposition ~\ref{projicerat} gives a complete description of $R$ except for the functions $C_{(K,\tau)}$. As in the Artinian case we are particularly interested in whether the $C_{(K,\tau)}$ are zero or not. Indeed, if we denote by $M_{(K,\tau)}$ the ideal generated by 
$\{z_i^{\alpha_i};~~i\in K,~~\alpha_\tau=(\alpha_1, \ldots, \alpha_n)\}$,  
then  
\[
\ann \dbar \Big[\frac{1}{z_{k_1}^{\alpha_{k_1}}}\Big ]\wedge\ldots\wedge\dbar \Big [\frac{1}{z_{k_\ell}^{\alpha_{k_\ell}}}\Big ] = M_{(K,\tau)},
\]
and so $\ann R_{(K,\tau)}=M_{(K,\tau)}$ if $C_{(K,\tau)}\not \equiv 0$. Note that $\beta$ does not affect the annihilator. Furthermore, to annihilate $R^{\p_K}$ one has to annihilate each entry $R_{(K,\tau)}$ and therefore 
\[
\ann R^{\p_K}=\bigcap _{(K,\tau);~~ C_{(K,\tau)}\not \equiv 0} M_{(K,\tau)}.
\]
Thus, in light of \eqref{bach}, a necessary condition for $C_{(K,\tau)}$ to be not identically zero is that $M\subset M_{(K,\tau)} $ and as in the Artinian case it turns out that we can determine precisely for which pairs $(K,\tau)$ this happens if $M$ is generic.

To this end we need to recall from ~\cite{BPS} how one can find the irredundant irreducible decomposition of a generic monomial ideal $M$ from the Scarf complex of a certain associated ideal. When $M$ is Artinian all irreducible components in the decomposition are of course $\m$-primary and we saw above, in Section  ~\ref{artinian_case}, that they correspond to facets of the Scarf complex $\Delta_M$. The idea is now that adding to $M$ a set of ``ghost generators'' $\{z_i^D\}$, where $D$ is some integer larger than the degree of any generator $m_i$, enables us to identify a $\p_K$-primary component $(z_{k}^{\alpha_k})_{k\in K}$ in the irreducible decomposition of $M$ with the $\m$-primary irreducible ideal generated by $\{z_{k}^{\alpha_k}\}_{k\in K}\cup \{z_i^D\}_{i\notin K}$ which in turn corresponds to a facet of the Scarf complex of the Artinian ideal $M^*:=(m_1,\ldots, m_r, z_1^D, \ldots, z_n^D)$. More precisely, for each subset $J$ of the underlying vertex set $\{1,\ldots, r,1_{\text{ghost}}, \ldots, n_{\text{ghost}}\}$ with labels $\{m_1,\ldots, m_r, z_1^D, \ldots, z_n^D\}$ let $M_J$ be the irreducible ideal generated by 
$\{z_i^{\alpha_i}, \alpha_i=\deg_{z_i}(m_J), \alpha_i<D\}$. 
Now Theorem 3.7 in ~\cite{BPS} states that the irredundant irreducible decomposition of a generic monomial ideal $M$ is given as the intersection of the irreducible ideals $M_J$, where $J$ runs over all facets of $\Delta_{M^*}$. To be able to use this result let us observe that $J\subset \{1,\ldots, r,1_{\text{ghost}}, \ldots, n_{\text{ghost}}\}$ can be identified with a pair $(K,\tau)$ above, by letting $K$ be determined by the set of ghost vertices and $\tau$ by the remaining vertices.  Indeed, given $J$, let $K=\{i\in\{1,\ldots,n\}; i_{\text{ghost}}\notin J\}$ and $\tau=\{i\in\{1,\ldots,r\}; i\in J\}$. With this identification the ideals $M_J$ and $M_{(K,\tau)}$ coincide.

We have the following generalization of Proposition ~\ref{coefficients}. 
\begin{prop}\label{lescarf}
Let $M\subset S$ be a generic monomial ideal and let $R$ be the residue current associated with the polyhedral cell complex $X$. Suppose that $\mathbb F_X$ is a resolution of $S/M$. Then $C_{(K,\tau)} \not \equiv 0$ if and only if $(K,\tau)$ is a facet of the Scarf complex $\Delta_{M^*}$.
\end{prop}
\begin{proof}
Observe that $M^*$ is generic. Thus, if $(K,\tau)$ is not a facet of $\Delta_{M^*}$, then there is a (minimal) generator $m'$ of $M^*$ such that $m'$ strictly divides $m_{(K,\tau)}$, so that $m'\notin \m^{\alpha_{(K,\tau)}}$, as we showed in the proof of Proposition ~\ref{coefficients}. In particular, $m'\notin M_{(K,\tau)}$, since clearly $M_{(K,\tau)}\subset \m^{\alpha_{(K,\tau)}}$. Moreover $m'$ has to be in $M$. Indeed, it is easy to see that $m'$ can not possibly be any of the generators $z_i^D$ of $M^*$. Since $D$ is chosen to be larger than the degree of any minimal generator of $M$, the degree of any of the variables in $m_{(K,\tau)}$ can not exceed $D$, and so $m'$ that strictly divides $m_{(K,\tau)}$ can not have degree $D$ in any variable.
Hence, if $(K,\tau)$ is not a facet of $\Delta_{M^*}$, then $M\not\subseteq M_{(K,\tau)}$ and so $C_{(K,\tau)}$ has to be 0.

On the other hand, since \eqref{irred} is irredundant, $C_{(K,\tau)}$ has to be not identically equal to zero whenever $(K,\tau)$ is a facet of $\Delta_{M^*}$. 
\end{proof}

To conclude, Propositions ~\ref{projicerat} and ~\ref{lescarf} give the following description of the residue current associated with a cellular resolution of a generic ideal, generalizing Theorem  ~\ref{main_artin}.
\begin{thm}\label{frost}
Let $M\subset S$ be a generic monomial ideal and let $R$ be the residue current associated with the polyhedral cell complex $X$. Suppose that $\mathbb F_X$ is a resolution of $S/M$. Then 

\begin{equation}\label{dickinson}
R=\sum_{(K,\tau) \text{ facet of } \Delta_{M^*}} 
R_{(K,\tau)} ~ e_\tau,
\end{equation}  
where $R_{(K,\tau)}$ is given by \eqref{himmelsgron} and $C_{(K,\tau)}$ there is not identically equal to zero.
\end{thm}
Note that \eqref{dickinson} is a refinement of the decomposition \eqref{pucko}, corresponding to that irreducible decompositions of monomial ideals are refinements of primary decompositions. 
In fact, Theorem ~\ref{dickinson} implies that the irredundant irreducible decomposition of a generic monomial ideal $M$ is recovered as 
\begin{equation}\label{emily}
M=\bigcap _{(K,\tau) \text{ facet of } \Delta_{M^*}} \ann R_{(K,\tau)};
\end{equation}
compare to Corollary \ref{dieu}. 

As in the Artinian case, we should remark, that even though we can not determine the set of non-vanishing entries of the residue current associated with an arbitrary cellular resolution, we can still estimate the number of them from below by the number or irreducible components of the corresponding ideal.

Finally, let us illustrate our results by some examples.

\begin{ex}\label{unkel}
Let $M$ be the (generic) ideal
\[
M=(z_1^4, z_1^2z_2, z_1z_2^2)=:(m_1, m_2, m_3),
\]
and let $M^*$ be the corresponding Artinian ideal obtained by adding the ghost vertices $z_1^D$ and $z_2^D$ (in fact, only $z_2^D$ comes into account), see Figure ~\ref{nonartinian}.

\begin{figure}\label{nonartinian}
\begin{center}
\includegraphics{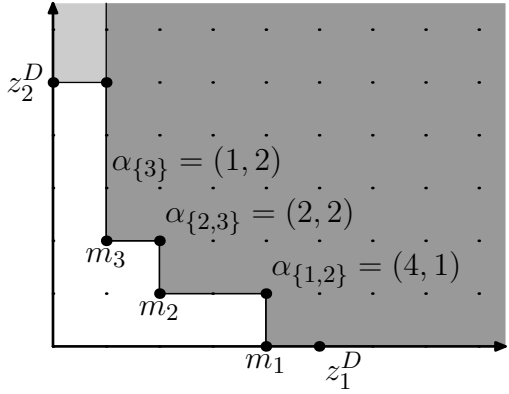}
\caption{The ideals $M$ (dark grey) and $M^*$ (light grey) in Example \ref{unkel}.}
\end{center}
\end{figure}
As in Example ~\ref{lexemple} the Scarf complex $\Delta_{M^*}$ is one-dimensional and the facets are the pairs of adjacent generators, that is 
$\tau_1=\{z_1^4, z_1^2z_2\}$, $\tau_2=\{z_1^2z_2, z_1z_2^2\}$ and $\tau_3=\{z_1z_2^2, z_2^D\}$. Note that the associated prime ideals of $M$ are $\p_{\{1\}}=(z_1)$ and $\p_{\{1,2\}}=(z_1,z_2)$, corresponding to $K=\{1\}$ and $\{1,2\}$ respectively.

Let $R$ be the residue current obtained from a cellular resolution. By Theorem ~\ref{frost} the current $R^{\p_{\{1\}}}$ is equal to a sum of terms
$R_{(\{1\},\tau)}~e_\tau$, 
where the sum is taken over facets $(\{1\},\tau)$ of $\Delta_{M^*}$. Recall that $(K,\tau)$ should be interpreted as the facet of $\Delta_{M^*}$ that has the vertex set $\{z_j^D\}_{j\notin K }\cup \{m_i\}_{i\in \tau}$, and thus  
we are looking for facets containing the ghost vertex $z_2^D$. However, there is only one such facet, namely $\tau_3=(\{1\},\{3\})$. Moreover, $\alpha_{\{3\}}=(1,2)$, and so  
\[
R_{\tau_3}=(\dbar \sigma_2) C_3(z_2) \otimes \dbar\bigg [\frac{1}{z_1} \bigg],
\]
with annihilator equal to $(z_1)$, which is precisely $M_{\tau_3}$.

Next, $R^{\p_{\{1,2\}}}=\sum R_{(\{1,2\},\tau)} ~e_\tau$, 
where the sum now is taken over facets of $\Delta_{M^*}$ that contain no ghost vertices. There are two such facets, $\tau_1=(\{1,2\},\{1,2\})$ and $\tau_2=(\{1,2\},\{2,3\})$, with corresponding currents 
\[
R_{\tau_1}=C_1 ~\dbar\bigg [\frac{1}{z_1^4} \bigg]\wedge\dbar\bigg [\frac{1}{z_2} \bigg]
\text{ and }
R_{\tau_2}=C_2 ~\dbar\bigg [\frac{1}{z_1^2} \bigg]\wedge\dbar\bigg [\frac{1}{z_2^2} \bigg].
\]
The annihilators are $(z_1^4,z_2)=M_{\tau_1}$ and 
$(z_1^2,z_2^2)=M_{\tau_2}$, respectively. Hence, according to \eqref{emily}, the 
irredundant irreducible decomposition is given by
\[
M=(z_1)\cap(z_1^4,z_2)\cap (z_1^2,z_2^2), 
\]
which could of course be seen directly.
\end{ex}

\begin{ex}
In ~\cite{AW} residue currents were used to obtain the following version of the Ehrenpreis-Palamodov fundamental principle: any smooth solution to the system of equations 
\begin{equation}\label{diff}
\eta (i \partial/\partial t) \cdot \xi(t)=0, \eta\in J\subset S
\end{equation}
on a smoothly bounded convex set in $\R^n$ can be written 
\begin{equation*}\label{losning}
\xi(t)=\int_{\C^n}\sum_k R_\ell^T(z) A_\ell(z) e^{-i\langle t,z\rangle},
\end{equation*}
for appropriate explicitly given vectors $A_\ell$ of smooth functions. Here $R_\ell^T$ are the (the transpose of) the components of the residue current associated with $J$. Conversely, any $\xi(t)$ given in this way is a homogeneous solution since $J=\ann R$.

Suppose that $J$ is a generic monomial ideal. Then, according to Theorem ~\ref{frost}, the solutions to \eqref{diff} are of the form 
\begin{equation}\label{diffa}
\xi(t)=\sum_{(K,\tau) \text{ facet of } \Delta_{M^*}} 
\int_{\C^n} 
\dbar \Big[\frac{1}{z_{k_1}^{\alpha_{k_1}}}\Big ]\wedge\ldots\wedge\dbar \Big[\frac{1}{z_{k_\ell}^{\alpha_{k_\ell}}}\Big ]\wedge
 A_{(K,\tau)}(z) e^{-i\langle t,z\rangle},
\end{equation}
where $A_{(K,\tau)}$ is smooth outside a set of codimension $>|K|$ and $(\alpha_1,\ldots,\alpha_n)=\alpha_\tau$. 
Now, the term in \eqref{diffa} indexed by $(K,\tau)$ is a polynomial in $\{z_i\}_{i\in K}$ of multidegree strictly smaller than $\alpha_\tau$ times a quite arbitrary function in $\{z_i\}_{i\notin K}$. 
It is easily checked directly that the general solution is a superposition of such functions, compare to for example Chapter 10.4 in ~\cite{St}.  
\end{ex}

{\bf Acknowledgements:} I would like to thank Mats Andersson for interesting discussions on the topic of this paper and for valuable comments on preliminary versions. Thanks to Ezra Miller for illuminating discussions on these matters and for the help with finding a result I needed. Thanks also to the referee for many helpful suggestions.

\def\listing#1#2#3{{\sc #1}:\ {\it #2},\ #3.}

\end{document}